\documentclass{amsart} 
\usepackage{amssymb}
\usepackage{amsmath}
\usepackage{amsfonts}

\begin{document}
\newtheorem{theo}{Theorem}[section]
\newtheorem{prop}[theo]{Proposition}
\newtheorem{lemma}[theo]{Lemma}
\newtheorem{coro}[theo]{Corollary}
\theoremstyle{definition}
\newtheorem{exam}[theo]{Example}
\newtheorem{defi}[theo]{Definition}
\newtheorem{rem}[theo]{Remark}


\newcommand{\Bb}{{\bf B}}
\newcommand{\Cb}{{\bf C}}
\newcommand{\Nb}{{\bf N}}
\newcommand{\Qb}{{\bf Q}}
\newcommand{\Rb}{{\bf R}}
\newcommand{\Zb}{{\bf Z}}
\newcommand{\Ac}{{\mathcal A}}
\newcommand{\Bc}{{\mathcal B}}
\newcommand{\Cc}{{\mathcal C}}
\newcommand{\Dc}{{\mathcal D}}
\newcommand{\Fc}{{\mathcal F}}
\newcommand{\Ic}{{\mathcal I}}
\newcommand{\Jc}{{\mathcal J}}
\newcommand{\Kc}{{\mathcal K}}
\newcommand{\Lc}{{\mathcal L}}
\newcommand{\Mx}{{\mathcal M}}
\newcommand{\Nc}{{\mathcal N}}
\newcommand{\Oc}{{\mathcal O}}
\newcommand{\Pc}{{\mathcal P}}
\newcommand{\Qc}{{\mathcal Q}}
\newcommand{\Sc}{{\mathcal S}}
\newcommand{\Tc}{{\mathcal T}}
\newcommand{\Uc}{{\mathcal U}}
\newcommand{\Vc}{{\mathcal V}}
\newcommand{\btu}{\bigtriangleup}

\author{Charles Akemann and Nik Weaver}

\title [Pure states on $\Bc(H)$]
       {Not all pure states on $\Bc(H)$ are diagonalizable}

\address {Department of Mathematics\\
          University of California\\
          Santa Barbara, CA 93106}
\address {Department of Mathematics\\
          Washington University in Saint Louis\\
          Saint Louis, MO 63130}

\email {akemann@math.ucsb.edu, nweaver@math.wustl.edu}


\date{June 6, 2006}

\begin{abstract}
Assuming the continuum hypothesis, we prove that $\Bc(H)$ has
a pure state whose restriction to any masa is not pure.
This resolves negatively an old conjecture of Anderson.
\end{abstract}

\maketitle


Let $H$ be a separable infinite-dimensional Hilbert space and let
$\Bc(H)$ be the algebra of bounded operators on $H$.  Anderson \cite{A2}
conjectured that every pure state on $\Bc(H)$ is {\it diagonalizable},
i.e., of the form $f(A) = \lim_\Uc \langle Ae_n, e_n\rangle$ for some
orthonormal basis $(e_n)$ and some ultrafilter $\Uc$ over $\Nb$.

A {\it masa} of $\Bc(H)$ is a maximal abelian self-adjoint subalgebra, and
an {\it atomic masa} is the set of all operators which are diagonalized
with respect to some given orthonormal basis of $H$.
Anderson's conjecture is related to a fundamental problem in C*-algebra, the
Kadison-Singer problem \cite{KS}, which asks whether every pure state
on an atomic masa of $\Bc(H)$ has a unique extension to a pure state
on $\Bc(H)$. If $(e_n)$ is an orthonormal basis of $H$, then every
pure state $f_0$ on the corresponding atomic masa $\Mx$ has the form
$f_0(A) = \lim_\Uc \langle Ae_n, e_n\rangle$ for some ultrafilter $\Uc$
over $\Nb$ and all $A \in \Mx$, and Anderson \cite{A1} showed that the
same formula, now for $A \in \Bc(H)$, defines a pure state $f$ on $\Bc(H)$.
Thus, a positive solution to the Kadison-Singer problem would say that
$f$ is the only pure state on $\Bc(H)$ which extends $f_0$.

In the presence of a positive solution to the Kadison-Singer problem,
Anderson's conjecture is equivalent to the weaker statement that every
pure state on $\Bc(H)$ restricts to a pure state on some atomic masa.
However, assuming the continuum hypothesis, we show that this weaker
statement is false: in fact, there exist pure states on $\Bc(H)$ whose
restriction to any masa is not pure. It follows that there are
pure states on $\Bc(H)$ that are not diagonalizable. It seems likely
that the statement ``every pure state on $\Bc(H)$ restricts to a pure
state on some atomic masa'' is also consistent with standard set theory.
This together with a positive solution to the Kadison-Singer
problem would imply the consistency of a positive answer to
Anderson's conjecture.

The key lemma we need is the following. Let $\Kc(H)$ be
the algebra of compact operators on $H$, let $\Cc(H) = \Bc(H)/\Kc(H)$
be the Calkin algebra, and let $\pi: \Bc(H) \to \Cc(H)$ be the natural
quotient map. We also write $\dot{a}$ for $\pi(a)$, for any $a \in \Bc(H)$.

\begin{lemma}\label{LM}
Let $\Ac$ be a separable C*-subalgebra of $\Bc(H)$ which
contains $\Kc(H)$, let $f$ be a pure state on $\Ac$ that annihilates
$\Kc(H)$, and let $\Mx$ be a masa of $\Bc(H)$. Then there is a pure
state $g$ on $\Bc(H)$ that extends $f$ and whose restriction to
$\Mx$ is not pure.
\end{lemma}

\begin{proof}
By Proposition 6 of \cite{A0} we can find an infinite-rank projection
$p \in \Bc(H)$ such that
\begin{equation}
\dot{p}\dot{a}\dot{p} = f(a)\dot{p}
\label{AA}
\end{equation}
for all $a \in \Ac$.

Lemma 1.4 and Theorem 2.1 of \cite{JP} imply that
$\pi(\Mx)$ is a masa of $\Cc(H)$. It follows that there
is a projection $q \in \Mx$ such that $\dot{q}$ neither contains nor is
orthogonal to $\dot{p}$. Otherwise $\dot{p}$ would be in the commutant
of $\pi(\Mx)$, and hence would belong to $\pi(\Mx)$ by maximality. But
this would mean $\dot{p}$ is minimal in $\pi(\Mx)$ because
any nonzero projection below $\dot{p}$ neither contains nor is orthogonal
to $\dot{p}$, and $\pi(\Mx)$ has no minimal projections.

Let $\phi: \Cc(H) \to \Bc(K)$ be an irreducible representation of the
Calkin algebra. It
is faithful because $\Cc(H)$ is simple. Therefore $\phi(\dot{q})$ neither
contains nor is orthogonal to $\phi(\dot{p})$, so we can find a unit
vector $v \in K$ in the range of $\phi(\dot{p})$ which is neither
contained in nor orthogonal to the range of $\phi(\dot{q})$. Finally,
define $g(a) = \langle \phi(\dot{a})v, v\rangle$ for all $a \in \Bc(H)$.
This is a pure state on $\Bc(H)$ because $\phi \circ\pi$ is an irreducible
representation of $\Bc(H)$. It extends $f$ because, using (\ref{AA}),
$$g(a) = \langle \phi(\dot{a})v, v\rangle
= \langle \phi(\dot{a})\phi(\dot{p})v, \phi(\dot{p})v\rangle
= \langle \phi(\dot{p}\dot{a}\dot{p})v, v\rangle
= \langle f(a)\phi(\dot{p})v, v\rangle
= f(a)$$
for all $a \in \Ac$. Finally, its restriction to $\Mx$ is not pure
because the projection $q \in \Mx$ has the property that
$$g(q) = \langle \phi(\dot{q})v, v\rangle$$
is strictly between $0$ and $1$, since $v$ is neither contained in
nor orthogonal to the range of $\phi(\dot{q})$.
\end{proof}

\begin{theo}\label{TH}
Assume the continuum hypothesis. Then there is a pure state on
$\Bc(H)$ whose restriction to any masa is not pure.
\end{theo}

\begin{proof}
Let $(a_\alpha)$, $\alpha < \aleph_1$, enumerate the elements of $\Bc(H)$.
Since every von Neumann subalgebra of $\Bc(H)$ is countably generated, a
simple cardinality argument shows that there are only $\aleph_1$ such
subalgebras. Hence $\Bc(H)$ has only $\aleph_1$ masas.
Let $(\Mx_\alpha)$, $\alpha < \aleph_1$, enumerate the masas of $\Bc(H)$.

We now inductively construct a nested transfinite sequence of unital
separable C*-subalgebras $\Ac_\alpha$ of $\Bc(H)$ together with pure
states $f_\alpha$ on $\Ac_\alpha$ such that for all $\alpha < \aleph_1$
\begin{enumerate}
\item
$a_\alpha \in \Ac_{\alpha+1}$
\item
if $\beta < \alpha$ then $f_\alpha$ restricted to $\Ac_\beta$ equals
$f_\beta$
\item
$\Ac_{\alpha+1}$ contains a projection $q_\alpha \in \Mx_\alpha$ such that
$0 < f_{\alpha+1}(q_\alpha) < 1$.
\end{enumerate}
Begin by letting $\Ac_0$ be any separable C*-subalgebra of $\Bc(H)$
that is unital and contains $\Kc(H)$ and let $f_0$ be any pure state on
$\Ac_0$ that annihilates $\Kc(H)$. At successor stages, use the lemma
to find a projection $q_\alpha \in \Mx_\alpha$ and a pure state
$g$ on $\Bc(H)$ such that $g|_{\Ac_\alpha} = f_\alpha$ and
$0 < g(q_\alpha) < 1$. By (\cite{AW}, Lemma 4) there is a separable
C*-algebra $\Ac_{\alpha+1} \subseteq \Bc(H)$ which contains $\Ac_\alpha$,
$a_\alpha$, and $q_\alpha$ and such that the restriction $f_{\alpha+1}$
of $g$ to $\Ac_{\alpha+1}$ is pure. Thus the construction may proceed.
At limit ordinals $\alpha$, let $\Ac_\alpha$ be the closure of
$\bigcup_{\beta<\alpha} \Ac_\beta$. The state $f_\alpha$ is determined
by the condition $f_\alpha|_{\Ac_\beta} = f_\beta$, and it is easy to
see that $f_\alpha$ must be pure. (If $g_1$ and $g_2$ are states on
$\Ac_\alpha$ such that $f_\alpha = (g_1 + g_2)/2$, then for all
$\beta < \alpha$ purity of $f_\beta$ implies that $g_1$ and $g_2$
agree when restricted to $\Ac_\beta$; thus $g_1 = g_2$.) This
completes the description of the construction.

Now define a state $f$ on $\Bc(H)$ by letting $f|_{\Ac_\alpha} = f_\alpha$.
By the reasoning used immediately above, $f$ is pure, and since
$0 < f(q_\alpha) < 1$ for all $\alpha$, the restriction of $f$
to any masa is not pure.
\end{proof}

It is interesting to contrast Theorem \ref{TH} with Theorem 9 of
\cite{A0}, which states that (assuming the continuum hypothesis)
any state on $\Cc(H)$ restricts to a pure state on some masa of
$\Cc(H)$. This does not conflict with our result because there are
many masas of $\Cc(H)$ which do not come from masas of $\Bc(H)$
(regardless of the truth of the continuum hypothesis). Indeed, $\Bc(H)$
has $2^{\aleph_0}$ masas but $\Cc(H)$ has $2^{2^{\aleph_0}}$ masas.
This can be seen by first finding $2^{\aleph_0}$ mutually orthogonal
nonzero projections $p_\alpha$ in $\Cc(H)$ \cite{W}, then finding
projections $q^1_\alpha, q^2_\alpha < p_\alpha$
such that $q^1_\alpha q^2_\alpha \neq q^2_\alpha q^1_\alpha$ for each
$\alpha$, and finally for each set $S \subseteq 2^{\aleph_0}$ choosing a
masa of $\Cc(H)$ that contains $\{q^1_\alpha: \alpha \in S\}$ and
$\{q^2_\alpha: \alpha \not\in S\}$. It is easy to see that this produces
$2^{2^{\aleph_0}}$ distinct masas.



\bigskip
\bigskip
\begin{thebibliography}{aaaaaaaa}

\bibitem{AW}
C.\ Akemann and N.\ Weaver, Consistency of a counterexample to
Naimark's problem, {\it Proc.\ Nat.\ Acad.\ Sci.\ USA \bf 101}
(2004), 7522-7525.

\bibitem{A0}
J.\ Anderson, Pathology in the Calkin Algebra, {\it J.\ Operator
Theory \bf 2} (1979), 159-167.

\bibitem{A1}
{---------}, Extreme points in sets of positive linear maps on
$\Bc(H)$, {\it J.\ Funct.\ Anal.\ \bf 31} (1979), 195-217.

\bibitem{A2}
{---------}, A conjecture concerning the pure states of $\Bc(H)$
and a related theorem, in {\it Topics in Modern Operator Theory},
pp.\ 27-43, Birkha\"user, 1981.

\bibitem{JP}
B.\ E.\ Johnson and S.\ K.\ Parrott, Operators commuting with a
von Neumann algebra modulo the set of compact operators,
{\it J.\ Funct.\ Anal.\ \bf 11} (1972), 39-61.

\bibitem{KS}
R.\ V.\ Kadison and I.\ M.\ Singer, Extensions of pure states,
{\it Amer.\ J.\ Math.\ \bf 81} (1959), 383-400.

\bibitem{W}
E.\ Wofsey, $P(\omega)/$fin and projections in the Calkin algebra,
manuscript.

\end{thebibliography}
\end{document}